\documentclass[journal]{IEEEtran}
\usepackage[english]{babel}
\usepackage{amsmath}
\usepackage{amssymb}
\usepackage{amsthm}
\usepackage{color}

\newtheorem{thm}{Theorem}

\usepackage{graphicx}
\graphicspath{ {pics/} }

\allowdisplaybreaks
\usepackage{etoolbox}
\makeatletter
\patchcmd{\@makecaption}
{\scshape}
{}
{}
{}
\makeatletter
\patchcmd{\@makecaption}
{\\}
{.\ }
{}
{}
\makeatother

\makeatletter
\long\def\@makecaption#1#2{\ifx\@captype\@IEEEtablestring%
	\footnotesize\begin{center}{\normalfont\footnotesize #1}\\
		{\normalfont\footnotesize\scshape #2}\end{center}%
	\@IEEEtablecaptionsepspace
	\else
	\@IEEEfigurecaptionsepspace
	\setbox\@tempboxa\hbox{\normalfont\footnotesize {#1.}~~ #2}%
	\ifdim \wd\@tempboxa >\hsize%
	\setbox\@tempboxa\hbox{\normalfont\footnotesize {#1.}~~ }%
	\parbox[t]{\hsize}{\normalfont\footnotesize \noindent\unhbox\@tempboxa#2}%
	\else
	\hbox to\hsize{\normalfont\footnotesize\hfil\box\@tempboxa\hfil}\fi\fi}
\makeatother

\IEEEoverridecommandlockouts  
\begin{document}
	
	\title{PDE Traffic Observer Validated on Freeway Data}
	
	\author{Huan~Yu, Alexandre Bayen, Miroslav Krstic
	%,~\IEEEmembership{Member,~IEEE,} Mirolsav Krstic~\IEEEmembership{Fellow,~IEEE,}%
}

\author{Huan Yu$^{*}$, Qijian Gan$^{\dagger}$, Alexandre Bayen$^{\dagger}$, Miroslav Krstic$^{*}$
	\thanks{$^{*}$Huan Yu and Miroslav Krstic are with the Department of Mechanical and Aerospace Engineering,
		University of California, San Diego, 9500 Gilman Dr, La Jolla, CA 92093
		 {(email: huy015@ucsd.edu; krstic@ucsd.edu)}}
	\thanks{$^{\dagger}$Qijian Gan and Alexandre Bayen are with the Department of Civil and Environmental Engineering, University of California, Berkeley, CA 94720{(email: qgan@berkeley.edu; bayen@berkeley.edu)}}
}

\maketitle

\begin{abstract}    
	
		This paper develops boundary observer for estimation of congested freeway traffic states based on Aw-Rascle-Zhang (ARZ) partial differential equations (PDE) model. Traffic state estimation refers to acquisition of traffic state information from partially observed traffic data. This problem is relevant for freeway due to its limited accessibility to real-time traffic information. We propose a model-driven approach in which estimation of aggregated traffic states in a freeway segment are obtained simply from boundary measurement of flow and velocity without knowledge of the initial states. The macroscopic traffic dynamics is represented by the ARZ model, consisting of $2 \times 2$ coupled nonlinear hyperbolic PDEs for traffic density and velocity. Analysis of the linearized ARZ model leads to the study of a hetero-directional hyperbolic PDE model for congested traffic regime. Using spatial transformation and PDE backstepping method, we construct a boundary observer consisting of a copy of the nonlinear plant with output injections from boundary measurement errors. The output injection gains are designed for the estimation error system so that the exponential stability of the error system in the $L^2$ norm and finite-time convergence to zero are guaranteed. Numerical simulations are conducted to validate the boundary observer design for estimation of the nonlinear ARZ model. In data validation, we calibrate model parameters of the ARZ model and then use vehicle trajectory data to test the performance of the observer design. 
\end{abstract}

%===============================================================================

\section{Introduction}
 
 Traffic state estimation plays an important role in traffic management. In order to mitigate freeway traffic congestion, various control algorithms~\cite{Carlson}~\cite{KM1}~\cite{Papageorgiou}~\cite{Huan1}~\cite{Huan:full}~\cite{Huan:adapt}~\cite{zpri19}~\cite{zhang16} are developed for ramp metering or variable speed limit. However, their performance heavily relies on accurate measurement of traffic states on mainline freeways. Due to financial and technical limitations, it is difficult to measure traffic states on mainline freeways everywhere at all times. Therefore, it is important to estimate traffic states at places where detection is missing.

 The topic of traffic state estimation refers to foreseeing traffc states with partially observed traffic data and some prior knowledge of traffic.  Such a topic has been extensively studied and attracted a lot of attentions in recent decades. According to the comprehensive review in \cite{bayen}, approaches on traffic estimation fall into the following three categories: model driven, data driven, and streaming data driven. Among them, the model driven approach is the most popular one and has been widely used to solve various traffic estimation problems. As a first step in the model driven approach, traffic flow models are often used to describe traffic dynamics and are calibrated with historical data. Then state estimates are obtained based on the calibrated model and real-time data inputs. Therefore, it is crucial for traffic estimation to have a physical model that can describe freeway traffic dynamics accurately.
 
Freeway traffic dynamics in spatial and temporal domains are usually described using macroscopic models with aggregated variables of traffic density, velocity and flux. These aggregated variables average out small-scale noises of freeway traffic and can be directly measured by stationary/point-based sensors like loop detectors. Among the macroscopic models, the Lighthill-Whitham-Richards (LWR) model by \cite{LW} and \cite{R} is one of the most applied models. This model is a first-order scalar hyperbolic PDE of density, and can predict the propagation and dissipation of traffic shockwaves and represent fundamental phenomena of free and congested regime of traffic. Several studies in \cite{Cl}~\cite{Coi02}~\cite{Coi03}~\cite{kest} have used such a model for traffic states estimation due to its simplicity and efficiency in model calibration and numerical simulation. However, the LWR model fails to describe stop-and-go traffic, which is the oscillatory behavior of congested traffic. The main reason is because the static equilibrium density-velocity relation of the LWR model is unable to reproduce the non-equilibrium relation appearing in the stop-and-go traffic.

In order to address this limitation, second-order models are proposed to employ additional nonlinear hyperbolic PDE for traffic velocity, in addition to the density conservation equation. Therefore, deviations from the equilibrium traffic relation are allowed in the second-order model since dynamics of the velocity PDE is captured. The first well-known second-order model is the Payne-Whitham (PW) model by \cite{p}~\cite{w}. But it predicts negative traffic velocity and information propagation faster than traffic which is physically unrealistic. Later in \cite{AW} and \cite{Zhang}, the Aw-Rascle-Zhang model was proposed which successfully addresses the anisotropic behavior of traffic and corrects the PW model’s prediction of traffic waves. For this reason, ARZ model has been studied intensively for the stop-and-go traffic over the recent years ~\cite{Bayen}~\cite{Seibold3}~\cite{Seibold1}~\cite{Helbing}~\cite{kerner}~\cite{Seibold2}. 
 
In the literature, there have been studies applying second-order models as physical models for traffic state estimation, for example, the second-order extended cell-transmission model in \cite{Mihaylova} and the second-order PW model in~\cite{Wa}. However, to the best of our knowledge, the nonlinear ARZ PDE model has never been used for state estimation. In order to accurately estimate the non-equilibrium traffic states for congested traffic, this paper applies the second-order ARZ model for observer design and data validation.

  In dealing with the second-order coupled nonlinear hyperbolic system, PDE control of the ARZ model has been studied through many recent efforts including~\cite{Bayen}~\cite{KBM}~\cite{KM1}~\cite{Huan1}~\cite{Huan2}~\cite{zpri}~\cite{zpri19}. The previous work by authors~\cite{Huan1}~\cite{Huan2} firstly consider adopts the PDE backstepping methodology for control of the ARZ model. Boundary control and observer design using PDE backstepping method have been developed for $2\times 2$ coupled hyperbolic systems~\cite{rafael}~\cite{Huan3} and the theoretical result for the general hetero-directional hyperbolic systems developed in~\cite{florent3}~\cite{florent1}~\cite{long}. The applications of the theoretical results include open-channel flow, oil drilling, heat exchangers and multi-phase flow problems, but have never been considered in traffic problems. 
  
  In \cite{Huan1}, an observer design is proposed for the linearized ARZ model in an effort to construct an output feedback controller. In \cite{Yu:observer}, we generalize the previous observer design to address the freeway traffic estimation problem from a more practical perspective. In specific, the observer design is proposed for the nonlinear ARZ model with certain assumptions of boundary conditions are removed. The observer design accepts a general functional form of the equilibrium density-velocity relation, rather than the Greenshield's model. The data validation results in this paper are obtained on the basis of the theoretical result in \cite{Yu:observer}. 
  
   In validation of the observer, vehicle trajectory data is used to obtain the aggregated values of traffic states. The ARZ model is calibrated with the historical field data. The model parameters are mostly obtained from historical data. The rest is determined from part of the dataset. Then the observer is constructed using the model parameters and real-time sensing of the data at boundaries. The performance of the PDE boundary observer is then evaluated with the field data in the temporal and spatial domain.

The contribution of this work: a systematic model-driven approach is developed for traffic state estimation. The PDE boundary observer based on the macroscopic ARZ traffic model is designed and validated. The theoretical observer design by backstepping method is generalized and adapted for the field-data validation. Vehicle trajectories data~\cite{NGSIM} is used to construct and to test the performance of the observer design. This result paves the way for implementing the PDE observer design in practice and give rise to a variety of opportunities to incorporate the PDE backstepping techniques in traffic estimation problem.

 The outline of this paper is as follows: 	
 in section II, we firstly introduce the nonlinear ARZ model, and analyze the linearized ARZ model for distinguishing the free and congested traffic. Section III designs the boundary observer for the linearized ARZ model using the backstepping method and the nonlinear boundary observer is developed using the output injections obtained from the linearized model. In section IV, numerical simulations of the nonlinear ARZ PDE model and state estimation by the nonlinear boundary observer are conducted firstly from an ad-hoc choice of model parameters. In section V, we calibrate the ARZ model with some field data and test the observer. The estimation errors are then analyzed.

\section{Problem statement}
We consider the traffic estimation problem for a stretch of freeway whose length is $L$. The macroscopic traffic dynamics is described by the ARZ model. We study the linearized ARZ model and discuss the characteristic speeds under the free and congested traffic regime. 
\subsection{Aw-Rascle-Zhang Model}
The ARZ model for $(x,t)\in [0,L]\times [{0, +\infty})$ is given 
\begin{align}
\partial_t \rho + \partial_x( \rho v)=&0,\label{rho} \\
\partial_t v+(v-\rho p'(\rho))\partial_x v=&\frac{V(\rho)-v}{\tau}. \label{v}
\end{align}
The state variable $\rho(x,t)$ denotes the traffic density and $v(x,t)$ denotes the traffic speed. The equilibrium velocity-density relationship $V(\rho)$ is a decreasing function of density. The equilibrium flux function $Q(\rho)$, also known as fundamental diagram, is defined as
\begin{align}
Q(\rho) = \rho V(\rho).
\end{align}
For ARZ model, the choice of $V(\rho)$ needs to satisfy that the flux function $Q(\rho)$ is smooth, strictly concave $Q(\rho)'' < 0$ and a strictly decreasing velocity functional form $V'(\rho) <0$. The second-order ARZ model is valid when the hyperbolicity is ensured for $Q(\rho)$. One of the basic choice of $V(\rho)$ is in the form of the Greenshield's model,
\begin{equation}
V(\rho)=v_f\left(1-\left(\frac{\rho}{\rho_m}\right)^\gamma\right). \label{vf}
\end{equation}
The observer design proposed in the following section is not limited by this choice. Later on, a more realistic functional form of $V(\rho)$ is proposed for the data validation which has a better fitting with traffic field data.

The inhomogeneous ARZ including a relaxation term on the right hand side of the velocity PDE is considered. The constant parameter $\tau$ is the relaxation time which describes drivers' driving behavior adapting to equilibrium density-velocity relation over time. Note that the homogeneous ARZ model without the relaxation term cannot address this phenomenon and poses an easier estimation problem. 

The increasing function of density $p(\rho)$ is defined as the traffic pressure
\begin{align}
p(\rho)=C_0\rho^\gamma,
\end{align} 
where $C_0,\gamma \in \mathbb{R}_+$, $p'(\rho)>0$ and $p(0)=0$ are assumed.
The pressure function $p(\rho)$ is chosen so that it is related to equilibrium velocity-density function $V(\rho)$ as 
\begin{align}
p(\rho) = V(0) - V(\rho).
\end{align}
Given $V(\rho)$ in \eqref{vf}, we have density pressure as
\begin{align}
p(\rho) = v_f\left( \frac{\rho}{\rho_m}\right)^\gamma.
\end{align}
Note that this choice of model parameter shows a marginal stability and a very slow damping effect of stop-and-go traffic. The following boundary observer design can be applied to the model when the above relation does not hold.

\subsection{Linearized ARZ model in traffic flux and velocity}
The traffic density is defined as the number of vehicles per unit length while the traffic flux represents the number of vehicles per unit time which cross a given point on the road. 
The traffic flow flux $q$ is defined as
\begin{align}
q=\rho v.
\end{align}
Traffic flux $q$ and velocity $v$ are the most accessible physical variables to measure in freeway traffic. $q$ is commonly measured by loop detectors and $v$ can be obtained by GPS or high-speed cameras. Therefore, we rewrite the ARZ model in traffic flux $q$ and traffic velocity $v$ as follows,
\begin{align}
\notag q_t + vq_x =&\frac{q}{v}\left(v+\frac{q}{v}V'\left(\frac{q}{v}\right)\right) v_x\\&+\frac{Q\left(\frac{q}{v}\right)-q}{\tau}, \label{nqv1}\\
v_t+\left(v+\frac{q}{v}V'\left(\frac{q}{v}\right)\right) v_x=&\frac{V\left(\frac{q}{v}\right)-v}{\tau}, \label{nqv2}
\end{align}
There is no explicit solution to the above nonlinear coupled hyperbolic system. To further understand the dynamics of the ARZ traffic model in $(q,v)$-system, we linearize the model around steady states ($q^\star$, $v^\star$) which are chosen as spatial and temporal nominal values of state variables. Small deviations from the nominal profile are defined as
\begin{align}
\tilde q(x,t)=&q(x,t)-q^\star,\\
\tilde v(x,t)=&v(x,t)-v^\star.
\end{align}
The steady density is given as $\rho^\star =q^\star/v^\star $ and setpoint density-velocity relation satisfy the equilibrium relation $V(\rho)$,
\begin{align}
v^\star = V(\rho^\star). 
\end{align}
The linearized ARZ model in $(\tilde q, \tilde v)$ around the reference system $(q^\star, v^\star)$ with boundary conditions is given by
	\begin{align}
	\notag\tilde{q}_t+\lambda_1\tilde{q}_x=&	\notag-\frac{q^\star}{v^\star}\left(v^\star + \frac{q^\star}{v^\star}V'\left(\frac{q^\star}{v^\star}\right)\right) \tilde{v}_x \\&-{q^\star}\frac{(v^\star)^2+q^\star V'\left(\frac{q^\star}{v^\star}\right)}{\tau (v^\star)^3}\tilde v +\frac{q^\star V'\left(\frac{q^\star}{v^\star}\right)}{\tau (v^\star)^2}\tilde q, \label{tq} \\
 \tilde v_t +  \lambda_2\tilde v_x =& -\frac{(v^\star)^2+q^\star V'\left(\frac{q^\star}{v^\star}\right)}{\tau (v^\star)^2}{\tilde v}+\frac{V'\left(\frac{q^\star}{v^\star}\right)}{\tau v^\star}\tilde q, \label{tv}
	\end{align}	
where the two characteristic speeds of the above linearized PDE model are
\begin{align}
\lambda_1 = &v^\star, \\ \lambda_2 =& v^\star + \frac{q^\star}{v^\star}V'\left(\frac{q^\star}{v^\star}\right) .
\end{align}
\begin{itemize}
	\item Free-flow regime  $: \lambda_1> 0  , \; \lambda_2 >0$\\
	In the free-flow regime, both the disturbances of traffic flux and velocity travel downstream, at respective characteristic speeds $\lambda_1$ and $\lambda_2 $. The linearized ARZ model in free-regime is a homo-directional hyperbolic system.
	\\
	\item Congested regime $: \lambda_1> 0  , \; \lambda_2 <0$\\ In the congested regime, the traffic density is greater than the critical value $\rho_c$ that satisfies $Q(\rho)'|_{\rho_c} = 0$ and the second characteristic speed $\lambda_2$ becomes the negative value. Therefore, disturbances of the traffic speed travel upstream with $\lambda_2$  while the disturbances of the traffic flow flux are carried downstream with the characteristic speed $\lambda_1$. The hetero-directional propagations of disturbances force vehicles into the stop-and-go traffic. 
\end{itemize}

In the free-flow regime, the linearized homo-directional hyperbolic PDEs can be solved explicitly by the inlet boundary values and therefore state estimates can be obtained by solving the hyperbolic PDEs. In this work, we focus on the congested regime with two hetero-directional hyperbolic PDEs. It is a more relevant and challenging problem for traffic states estimation.

\section{Boundary Observer Design}
In this section, boundary sensing is employed for the observer design. The state estimation of the nonlinear ARZ model is achieved using backstepping method. The output injection gains are designed for the linearized ARZ model and then adding to a copy of the nonlinear plant.

Boundary values of state variations from the steady states are defined as 
\begin{align}
Y_{q,in}(t)  = & \tilde q(0,t), \label{y1}\\
Y_{q,out}(t)  = & \tilde q(L,t),\\
Y_v(t)  = & \tilde v(L,t).\label{y2}
\end{align}
where the values of $\tilde q(0,t)$, $\tilde q(L,t)$ and $\tilde v(L,t)$ are obtained by subtracting setpoint values $(q^\star, v^\star)$ from the sensing of incoming traffic flux $q(0,t)$, outgoing flux $q(L,t)$ and outgoing velocity $v(L,t)$,
\begin{align}
y_q(t) =& q(0,t), \\
y_{out}(t)= & q(L,t),\\
y_v(t) = & v(L,t).
\end{align}
Sensing of the aggregated values of the traffic flux and velocity can be obtained by high-speed camera or induction loop detectors. The induction loops are coils of wire embedded in the surface of the road to detect changes of inductance when vehicles pass. The high-speed cameras record the vehicle trajectories for a freeway segment.
%The boundary conditions are assumed to be
%\begin{align}
%\rho(0,t) v(0,t) =& q^\star,\label{bc1}\\
%\rho(L,t) = & \rho^\star.\label{bc2}
%\end{align}
%where $q^\star = \rho^\star v^\star, v^\star = V(\rho^\star) $, $\rho^\star$ and $v^\star$ are steady states

\subsection{Output injection for linearized ARZ model}
We diagonalize the linearized equations and therefore write $(\tilde q, \tilde v)$-system in the Riemann coordinates. The Riemann variables are defined as
\begin{align}
\xi_1 = &\frac{\rho^\star \lambda_2}{\lambda_1-\lambda_2}\tilde v + \tilde q, \\ \xi_2 =& \frac{q^\star}{\lambda_1-\lambda_2}\tilde v.
\end{align}
The inverse transformation is given by 
\begin{align}
\tilde v =& \frac{\lambda_{1} - \lambda_{2}}{q^\star} \xi_2,\\
\tilde q =& \xi_1 - \frac{\lambda_{2}}{\lambda_{1}}\xi_2.
\end{align}
The measurements are taken at boundaries lead to the following boundary conditions
\begin{align}
\xi_1(0,t) =& \frac{\lambda_2}{\lambda_1} \xi_2(0,t) + Y_q(t),\\
\xi_2(L,t) = &  \frac{q^\star}{\lambda_1-\lambda_2} Y_v(t).
\end{align}
Therefore the linearized ARZ model in Riemann coordinates is obtained
\begin{align}
\partial_t \xi_1 + \lambda_1 \partial_x \xi_1 =& -\frac{1}{\tau}\xi_1,\\
\partial_t \xi_2 + \lambda_2 \partial_x \xi_2 =& -\frac{1}{\tau}\xi_1,\\
\xi_1(0,t) =& \frac{\lambda_2}{\lambda_1} \xi_2(0,t),\\
\xi_2(L,t) = & \xi_1(L,t).
\end{align}
%The initial conditions are obtained according to \eqref{ic1} and \eqref{ic2},
%\begin{align}
%\xi_1(x,0) =& \frac{q^\star}{\lambda_1-\lambda_2}\tilde v(x,0) + v^\star \tilde \rho (x,0),\\
%\xi_2(x,t) =& \frac{q^\star}{\lambda_1-\lambda_2}\tilde v(x,0).
%\end{align}
%where the variations are chosen to be the same with that of the nonlinear ARZ,
%\begin{align}
%\tilde \rho (x,0) = & 0.1 \sin\left(\frac{3\pi x}{L}\right)\rho^\star,\\
%\tilde v(x,0) = & -0.1 \sin\left(\frac{3\pi x}{L}\right)v^\star.
%\end{align}
In order to diagonalize the right hand side to implement the backstepping method, we introduce a scaled state as follows:
\begin{align}
\bar w(x,t) = & \exp\left(\frac{ x }{\tau \lambda_1 } \right)  \xi_1(x,t), \label{tr1}\\
\bar v(x,t) = & \xi_2(x,t).
\end{align}

The $(\xi_1, \xi_2 )$-system is then transformed to a first-order $2 \times 2$ hyperbolic system
\begin{align}
\bar w_t(x,t)  + \lambda_1 \bar w_x(x,t) =& 0, \label{bw}\\
\bar v_t(x,t) + \lambda_2 \bar v_x(x,t) =& c(x) \bar w(x,t), \label{bv}\\
\bar w(0,t)=&\frac{\lambda_2}{\lambda_1} \bar v(0,t) + Y_{q,in}(t) ,\label{bbvw0}\\
\bar v(L,t)=& \frac{q^\star}{\lambda_1-\lambda_2} Y_v(t), \label{bbvw1}
\end{align}
where the spatially varying parameter $c(x)$ is defined as 
\begin{align}
c(x)=-\frac{1}{\tau}\exp\left(-\frac{ x }{\tau \lambda_{1} } \right),
\end{align}	
Parameter $c(x)$ is a strictly increasing function and bounded by
\begin{align}
-\frac{1}{\tau}\leq c(x)\leq -\frac{1}{\tau}\exp\left(-\frac{L}{\tau \lambda_{1} } \right) . \label{cbound}
\end{align}

Then we design a boundary observer for the linearized ARZ model to estimate $\bar w(x,t)$ and $\bar v(x,t)$ by constructing the following system
\begin{align}
\hat w_t(x,t)  + \lambda_1 \hat w_x(x,t) =& r(x)(\bar w(L,t)-\hat w(L,t))  ,\label{O1}\\
\notag \hat v_t(x,t) + \lambda_2 \hat v_x(x,t) =& c(x) \hat w(x,t)\\
&+ s(x)(\bar w(L,t)-\hat w(L,t))  , \\ 
\hat w(0,t)=&
\frac{\lambda_2}{\lambda_1} \hat v(0,t) + Y_{q,in}(t),\\
\hat v(L,t)
=&\frac{q^\star}{\lambda_1-\lambda_2} Y_v(t) , \label{2}
\end{align}
where $\hat w(x,t)$ and $\hat v(x,t)$ are the estimates of the state variables $\bar w(x,t)$ and $\bar v(x,t)$. The value $\bar w(L,t)$ is obtained by plugging in the measured outgoing flow flux $Y_{q,out}(t)$ and velocity $Y_v(t)$ into \eqref{tr1},
\begin{align}
\bar w(L,t) = \exp\left(\frac{ L }{\tau \lambda_1 } \right)\left(\frac{\rho^\star \lambda_2}{\lambda_1-\lambda_2} Y_v(t) + Y_{q,out}(t)\right). \label{bwl}
\end{align}

The term $r(x)$ and $s(x)$ are output injection gains to be designed. We denote estimation errors as
\begin{align}
\check w(x,t)=&\bar w(x,t) -\hat w(x,t),\\
\check v(x,t)=&\bar v(x,t) -\hat v(x,t).
\end{align}
The error system is obtained by subtracting the estimates \eqref{O1}-\eqref{2} from \eqref{bw}-\eqref{bbvw1},
\begin{align}
\check w_t(x,t) + \lambda_1  \check w_x(x,t) =& r(x)\check w(L,t), \label{e21}\\
\notag	\check v_t(x,t) + \lambda_2  \check v_x(x,t)  =&c(x) \check w(x,t)+ s(x)\check w(L,t),\\ 
\check w(0,t)=&
\frac{\lambda_2}{\lambda_1} \check v(0,t), \\
\check v(L,t)
=& 0.\label{e24}
\end{align}

The design of output injection gains $r(x)$ and $s(x)$ needs to guarantee that the error system $(\check w, \check v)$ decays to zero. Using the backstepping transformation, we transform the error system \eqref{e21}-\eqref{e24} into the following target system
\begin{align}
\alpha_t(x,t) + \lambda_1 \alpha_x(x,t) =& 0, \label{al}\\
\beta_t(x,t) + \lambda_2 \beta_x(x,t) =& 0,\\ 
\alpha (0,t)=&\frac{\lambda_2}{\lambda_1} \beta(0,t),  \label{inb}\\
\beta(L,t)=& 0.\label{outb}
\end{align}
The explicit solution to the target system \eqref{al}-\eqref{outb} is easily found
\begin{align}
	\alpha(x,t) =& \alpha\left(0, t-\frac{x}{\lambda_1}\right), \quad\quad\;\; t> \frac{L}{|\lambda_{1}|},\\
	\beta(x,t) =& \beta \left(L, t+\frac{L-x}{\lambda_2}\right), \quad t> \frac{L}{|\lambda_{2}|}. 
\end{align}
Thus we have 
\begin{align}
\alpha(x,t)\equiv \beta(x,t)\equiv 0,
\end{align}
after a finite time $t=t_f$ where
\begin{align}
t_f = \frac{L}{|\lambda_{1}|} + \frac{L}{|\lambda_{2}|}. \label{tf}
\end{align}
It is straightforward to prove that the $\alpha, \beta$ system is $L^2$ exponentially stable. 

The backstepping transformation is given in the form of spatial Volterra integral
\begin{align}
\alpha(x,t) =& \check w(x,t)-\int^{L}_{x} K(L+x-\xi) \check w(\xi,t)d\xi, \label{bkst1}\\
\beta(x,t)=& \check v(x,t)-\int^{L}_{x} M(\lambda_1  x - \lambda_2\xi)\check w(\xi,t)d\xi,\label{bkst2}
\end{align}
where the kernel variables $K(x)$ and $M(x)$ map the error system into the target system where the coupling term on the right hand-side is eliminated by the output injections. The kernel $M(x)$ is defined as
\begin{align}
M(x)=-\frac{1}{\lambda_1  - \lambda_2} c\left(\frac{x}{\lambda_1 - \lambda_2}\right).
\end{align}
For boundary condition \eqref{inb} to hold, the kernels $ K(x)$ and $\check M(x)$ satisfy the relation
\begin{align}
K(L-\xi)=& M ((\lambda_2  - \lambda_1)\xi).
\end{align}
The kernel $K$ is then obtained
\begin{align}
K(x)=-\frac{1}{\lambda_1  - \lambda_2} c\left(\frac{-\lambda_2}{\lambda_1 - \lambda_2}(L-x)\right).
\end{align}
According to the boundedness of $c(x)$ in \eqref{cbound}, the kernels are bounded by
\begin{align}
|K(x)| \leq \frac{1}{(\lambda_1 - \lambda_2) \tau}, \label{bck}
\end{align}
and therefore $M(x)$ is bounded.
The output injection gain $r(x)$ and $s(x)$ are given by
\begin{align}
r(x)=&\lambda_1  K(x)=-\frac{\lambda_1}{\lambda_1 - \lambda_2} c\left(-\frac{ \lambda_2}{\lambda_1 - \lambda_2}(L-x)\right), \label{g1}\\
\notag	s(x)=&- \lambda_1 M( \lambda_1  x - \lambda_2 L)\\
=&\frac{\lambda_1}{\lambda_1 - \lambda_2}  c\left( x - \frac{ \lambda_2}{\lambda_1 - \lambda_2}(L-x)\right).\label{g2}
\end{align}
The backstepping transformation in \eqref{bkst1} and \eqref{bkst2} is invertible. Therefore, we study the stability of the error system through the target system \eqref{al}-\eqref{outb}. We arrive at the following theorem.

\begin{thm}
	Consider system \eqref{e21}-\eqref{e24} with inital conditions $\check w_0, \check v_0 \in L^2([0,L])$. The equilibrium $\check w \equiv \check v \equiv 0$ is exponentially stable in the $L^2$ sense. It holds that
	\begin{align}
	||\bar w(\cdot,t)-\hat w(\cdot, t)||\to 0\\
	||\bar v(\cdot,t)-\hat v(\cdot, t)||\to 0
	\end{align}
	and the convergence to the equilibrium is reached in the finite time $t=t_f$ given in \eqref{tf}.
\end{thm}

\subsection{Boundary observer design for Nonlinear ARZ model}
For nonlinear boundary observer, we construct the system by keeping the output injections that are designed for the linearized ARZ model, then add them to the copy of the original nonlinear ARZ model.

We summarize the transformation from the linearized ARZ model in $(\tilde q,\tilde v)$-system to $(\bar w, \bar v)$-system,
\begin{align}
\bar w(x,t)=& \exp\left(\frac{ x }{\tau \lambda_1 } \right)\left(\frac{\rho^\star \lambda_2}{\lambda_1-\lambda_2}\tilde v(x,t) + \tilde q(x,t)\right),\label{tf1}\\
\bar v(x,t)=&\frac{q^\star}{\lambda_1-\lambda_2}\tilde{v}(x,t).\label{tf2}
\end{align}
And the inverse transformation is given by
\begin{align}
\tilde q(x,t)=&\exp\left(-\frac{ x }{\tau  \lambda_1 } \right)\bar w(x,t) -  \frac{\lambda_2}{ \lambda_1}\bar v(x,t), \label{qv1}\\
\tilde v(x,t)=&\frac{\lambda_1-\lambda_2}{q^\star} \bar{v}(\xi,t).\label{qv2}
\end{align}
Due to the equivalence of $(\check w,\check v)$ and $(\tilde q, \tilde v)$-system, we arrive at the following theorem for the linearized ARZ model.
\begin{thm}
	Consider system \eqref{tq}-\eqref{tv} with inital conditions $\tilde q_0, \tilde v_0 \in L^2([0,L])$. The equilibrium $\tilde q \equiv \tilde v \equiv 0$ is exponentially stable in the $L^2$ sense. It holds that
	\begin{align}
	||q (\cdot,t)-q^\star||\to 0\\
	||v (\cdot,t)-v^\star||\to 0
	\end{align}
	and the convergence to set points is reached in finite time $t=t_f$.
\end{thm}

We denote the error injections designed for the linearized ARZ model \eqref{O1}-\eqref{2} as
\begin{align}
E_w(t) =& r(x)(\bar w(L,t)-\hat w(L,t)), \\
E_v(t) =& s(x)(\bar w(L,t)-\hat w(L,t)).
\end{align}
The output injection gains $r(x)$, $s(x)$ are designed in \eqref{g1} and \eqref{g2}. According to \eqref{bwl}, $\bar w(L,t)$ is obtained from the real-time measurement of the traffic boundary data in \eqref{y1}-\eqref{y2}. Therefore, the values of output injections $E_w(t)$ and $E_v(t)$ are known. 

The nonlinear observer for state estimation of density and velocity $(\hat\rho(x,t),\hat v(x,t))$ is obtained by combining the copy of the nonlinear ARZ model $(\rho,v)$ given by \eqref{rho}, \eqref{v} and the above linear injection errors in original state variables density and velocity,
\begin{align}
\partial_t \hat \rho + \partial_x( \hat \rho \hat v)=&  { \frac{1}{v^\star}\left(\exp\left(-\frac{ L }{\tau  \lambda_1 } \right)E_w- E_v\right)}, \\
\partial_t \hat v+(\hat v + \hat \rho V'(\hat \rho))\partial_x \hat v=&\frac{V(\hat \rho)- \hat v}{\tau} +{ \frac{\lambda_1-\lambda_2}{q^\star}  E_v},
\end{align}
where the linear injection on the right hand side are obtained by transforming $(\hat w, \hat v)$ to $(\rho, v)$ given in \eqref{qv1},\eqref{qv2}. The boundary conditions are
\begin{align}
\hat \rho(0,t) = \frac{y_q(t)}{\hat v(0,t)},\\
\hat v(L,t) = {y_v(t)}.
\end{align}
When the initial states of the system is close to the set points, the linearized part dominates the nonlinear estimation error system. Therefore $L^2$ exponential stability and finite-time convergence are achieved for the linearized ARZ model.
The local $H^2$ exponential stability can be derived for the estimation error system of the nonlinear ARZ model, following approach in \cite{rafael}. The estimation result is firstly validated in the following numerical simulation with an ad-hoc choice of parameters.

\begin{table}[t!]
	\caption{Parameter Table}
	\label{tab:default}
	\centering
	\begin{tabular}{ p{4cm}|p{2cm}  }	
		\hline
		%	\multicolumn{4}{|c|}{Country List} \\
		%	\hline
		\bf Parameter Name  & \bf Value \\
		\hline
		Maximum traffic density $\rho_m$   & 160 vehicles/km \\	\hline
		Traffic pressure and coefficient $\gamma$ & 1 \\	\hline
		Maximum traffic velocity $v_f$	&  40 m/s  \\	\hline
		Relaxation time $\tau$  & 60 s \\	\hline
		Reference density $\rho^\star$ & 120 vehicles/km \\	\hline
		Reference velocity $v^\star$ & 10 m/s\\	\hline
		Freeway segment length $L$ & $400$ m
		\\		\hline
	\end{tabular}
\end{table}

\section{Numerical Simulation}

For simulation of the nonlinear ARZ PDE model, we assume that
the initial conditions are sinusoidal oscillations around the steady states $(\rho^\star, v^\star)$ which are in the congested regime.
The initial conditions are assumed to be
\begin{align}
\rho(x,0) = &  0.1 \sin\left(\frac{3\pi x}{L}\right)\rho^\star + \rho^\star, \label{ic1}\\
v(x,0) = & -0.1 \sin\left(\frac{3\pi x}{L}\right) v^\star + v^\star. \label{ic2}
\end{align}
Model parameters of a one-lane traffic in the congested regime is considered and chosen as shown in the table 1. 

We consider a constant incoming flow and constant outgoing density for boundary conditions,
\begin{align}
\tilde q(0,t) = & 0, \label{bc1}\\
\tilde v(L,t) = & \frac{1}{\rho^\star}\tilde q(L,t).\label{bc2}	
\end{align}
In the next section, we validate the observer design with the traffic filed data, we do not prescribe any boundary conditions beforehand but directly take the measurement of the boundary data.

\begin{figure*}[t!]
	\includegraphics[width=13cm]{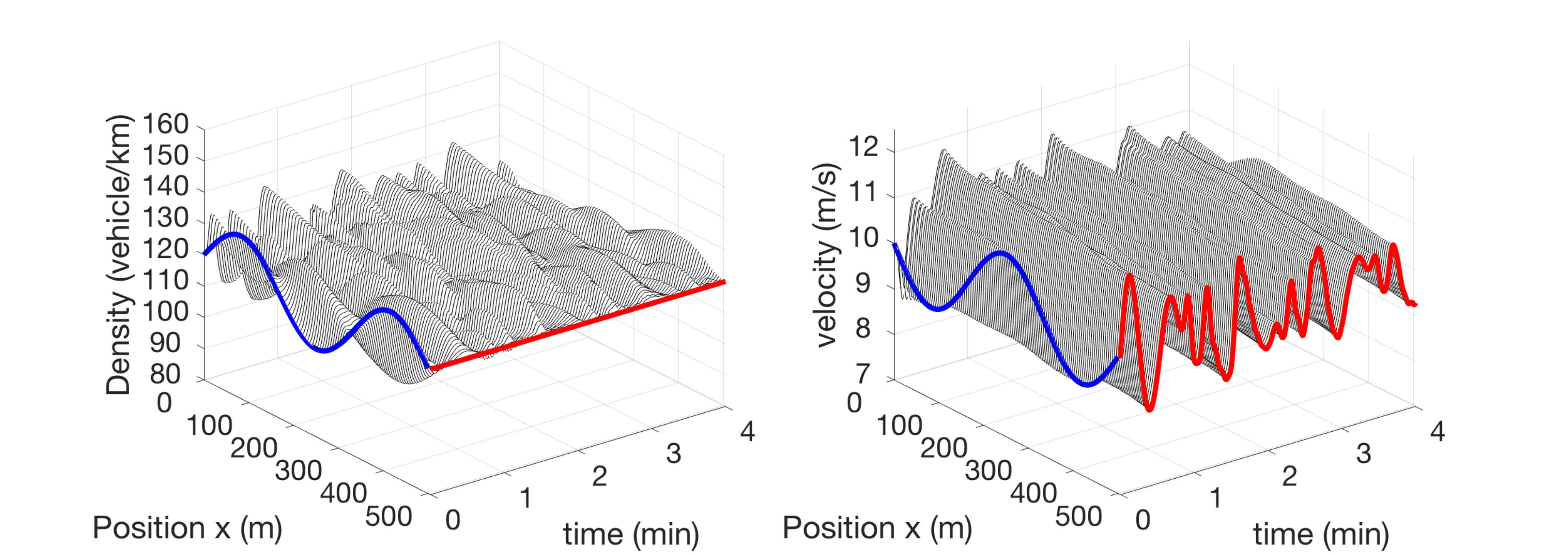}
	\centering
	\caption{Density $\rho(x,t)$ and velocity $v(x,t)$ of nonlinear ARZ model.}
\end{figure*} 
\begin{figure*}[t!]
	\includegraphics[width=13cm]{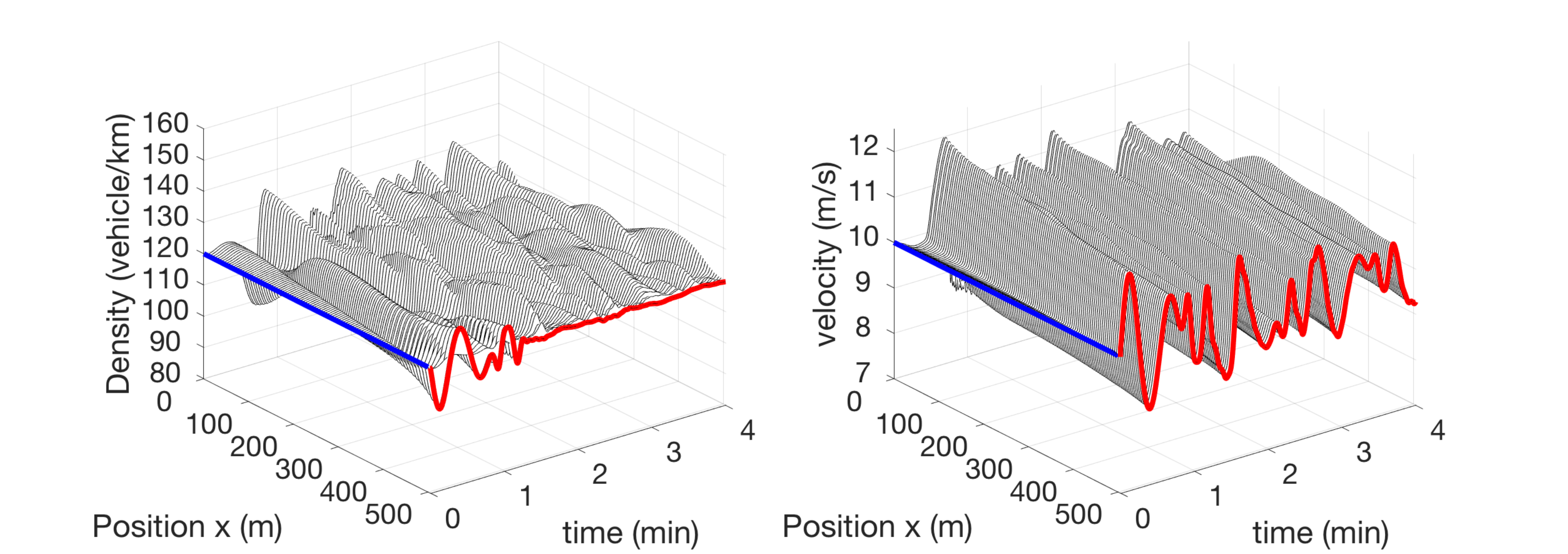}
	\centering
	\caption{States estimates $\hat \rho(x,t)$ and $\hat v(x,t)$ of nonlinear boundary observer.}
\end{figure*} 
\begin{figure*}[t!]
	\includegraphics[width=13cm]{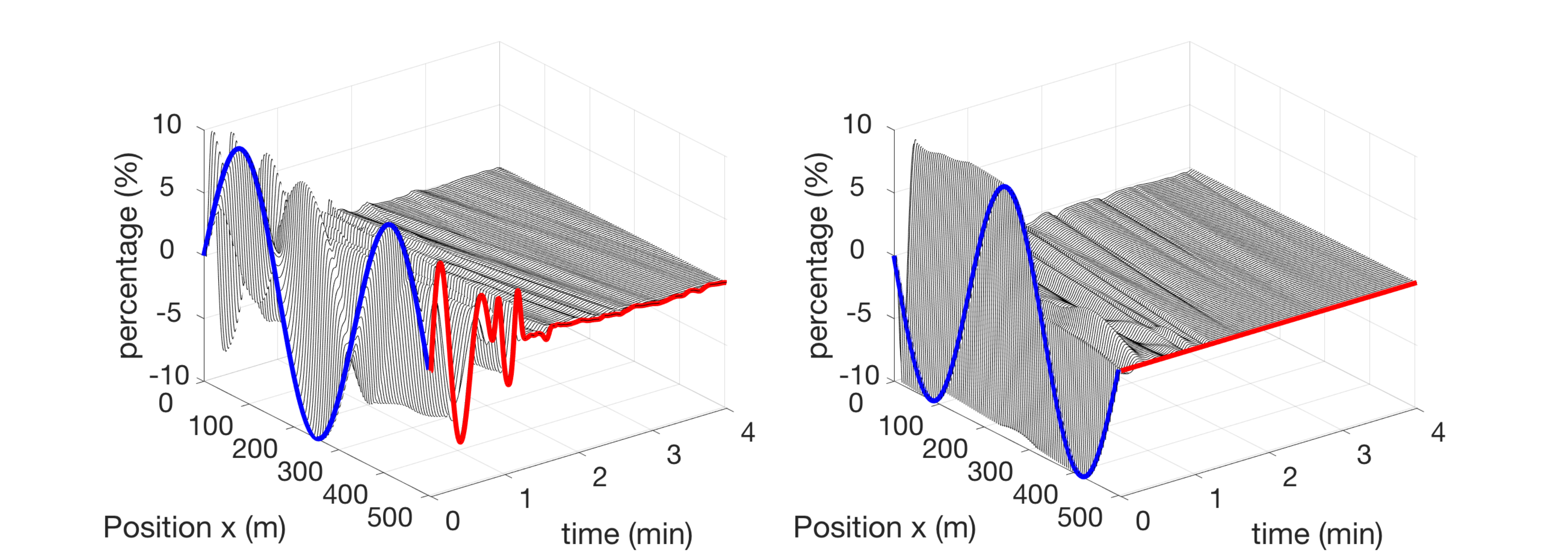}
	\centering
	\caption{Estimation errors $\check \rho(x,t)$ and $\check v(x,t)$.}
\end{figure*} 

We use the finite volume method which is common in traffic flow applications. The numerical approach divides the freeway segment into cells and then approximates the cell values considering the balance of fluxes through the boundaries of the adjacent cells. In order to obtain the numerical fluxes, we write the ARZ model in the conservative variables, then apply two-stage Lax-Wendroff scheme to discretize the ARZ model in spatio-temporal domain. The scheme is second-order accurate in space and first-order in time. The spatial grid resolution is chosen to be smaller than the average vehicle size so that the numerical errors are smaller than the model errors. Therefore the numerical simulation is valid for this continuum model.

The inhomogeneous nonlinear ARZ model written in the conservative form is given by
\begin{align}
\rho_t + (\rho v)_x =& 0,\\
y_t + (yv)_x =& -\frac{y}{\tau},
\end{align}
where $\rho$ and $y$ are conservative variables, and $y$ is defined as
\begin{align}
	y =  \rho(v - V(\rho)).
\end{align}
The numerical fluxes are then obtained by
\begin{align}
F_\rho =& y + \rho V(\rho),\\
F_y =& \frac{y^2}{\rho} + yV(\rho).
\end{align}
The Lax-wendroff numerical scheme is performed through two-stage update from $\left(\rho_j^n, y_j^n\right)$  to $\left(\rho_{j}^{n+1}, y_{j}^{n+1} \right)$.

At the first stage, the update law of $\left(\rho_j^n, y_j^n\right)$ to $\left(\rho_{j+\frac{1}{2}}^{n+\frac{1}{2}}, y_{j+\frac{1}{2}}^{n+\frac{1}{2}}\right)$ is given by
\begin{align}
\rho_{j+\frac{1}{2}}^{n+\frac{1}{2}} =& \frac{1}{2} \left(\rho_j^{n} + \rho_{j+1}^{n}\right) - \frac{\Delta t}{2\Delta x} \left((F_\rho)_{j+1}^{n}-(F_\rho)_{j}^{n}\right),\\
\notag y_{j+\frac{1}{2}}^{n+\frac{1}{2}}=& \frac{1}{2} \left(y_j^{n} + y_{j+1}^{n}\right) - \frac{\Delta t}{2\Delta x} \left((F_y)_{j+1}^{n}-(F_y)_{j}^{n}\right)\\& -\frac{ \Delta t}{4 \tau}\left(y_j^{n} + y_{j+1}^{n}\right),
\end{align}
Then we calculate the numerical flux at the intermediate points of state variables and the obtain the final stage as
\begin{align}
\rho_{j}^{n+1} =& \rho_{j}^{n} - \frac{\Delta t}{\Delta x}\left((F_\rho)_{j+\frac{1}{2}}^{n+\frac{1}{2}} -(F_\rho)_{j-\frac{1}{2}}^{n+\frac{1}{2}} \right),\\
\notag y_{j}^{n+1} =& y_{j}^{n} - \frac{\Delta t}{\Delta x}\left((F_y)_{j+\frac{1}{2}}^{n+\frac{1}{2}} -(F_y)_{j-\frac{1}{2}}^{n+\frac{1}{2}} \right) \\& -\frac{ \Delta t}{2 \tau}\left(y_{j+\frac{1}{2}}^{n+\frac{1}{2}} + y_{j-\frac{1}{2}}^{n+\frac{1}{2}}\right).
\end{align}
For the numerical stability of the Lax-Wendroff scheme, the spatial grid size $\Delta x$ and time step $\Delta t$ is chosen so that CFL condition is satisfied:
\begin{align}
\max |\lambda_{1,2}| \leq \frac{\Delta x}{\Delta t},
\end{align} 
%
%At the first stage, the update law of $\left(\rho_j^n, y_j^n\right)$ to $\left(\rho_{j+\frac{1}{2}}^{n+\frac{1}{2}}, y_{j+\frac{1}{2}}^{n+\frac{1}{2}}\right)$ is given 
%\begin{align}
%\rho_{j+\frac{1}{2}}^{n+\frac{1}{2}} =& \frac{1}{2} \left(\rho_j^{n} + \rho_{j+1}^{n}\right) - \frac{\Delta t}{2\Delta x} \left((F_\rho)_{j+1}^{n}-(F_\rho)_{j}^{n}\right),\\
%y_{j+\frac{1}{2}}^{n+\frac{1}{2}}=& \frac{1}{2} \left(y_j^{n} + y_{j+1}^{n}\right) - \frac{\Delta t}{2\Delta x} \left((F_y)_{j+1}^{n}-(F_y)_{j}^{n}\right) -\frac{ \Delta t}{4 \tau}\left(y_j^{n} + y_{j+1}^{n}\right),
%\end{align}
%Then we calculate the numerical flux at the intermediate points of state variables and the obtain the final stage as
%\begin{align}
%\rho_{j}^{n+1} =& \rho_{j}^{n} - \frac{\Delta t}{\Delta x}\left((F_\rho)_{j+\frac{1}{2}}^{n+\frac{1}{2}} -(F_\rho)_{j-\frac{1}{2}}^{n+\frac{1}{2}} \right),\\
%y_{j}^{n+1} =& y_{j}^{n} - \frac{\Delta t}{\Delta x}\left((F_y)_{j+\frac{1}{2}}^{n+\frac{1}{2}} -(F_y)_{j-\frac{1}{2}}^{n+\frac{1}{2}} \right) -\frac{ \Delta t}{2 \tau}\left(y_{j+\frac{1}{2}}^{n+\frac{1}{2}} + y_{j-\frac{1}{2}}^{n+\frac{1}{2}}\right).
%\end{align}

We specify state values at both $x=0$ and $x=L$ boundaries. ARZ model will pick up some combination of $\rho$ and $v$ at each of the two boundaries, depending on the direction of characteristics at the boundary cells. We implement the boundary conditions in \eqref{bc1} and \eqref{bc2}.

The numerical simulation result of the nonlinear ARZ, the nonlinear boundary observer estimation and the estimation errors are plotted in Fig. 1-3. Blue lines represent the initial conditions while the red lines represent the evolution of outlet state values in the temporal domain. The simulation is performed for a $500\;$m length of freeway segment and evolution of traffic states density and velocity are plotted for $4 \; \rm min$. 

In Fig. 1, traffic density and velocity are slightly damped and keeps oscillating in the domain. It takes the initial disturbance-generated vehicles to leave the domain in $50\;\rm s$ but the oscillations sustain for more than $4 \; \rm min$ which means the following incoming vehicles entering the acceleration-deceleration cycles under the influence of stop-and-go waves. The traffic states are chosen to be in the congested regime and the stop-and-go phenomenon is demonstrated in the simulation. 

State estimation of traffic density and velocity by the nonlinear observer is shown in Fig. 2. The measurement is taken for the outgoing velocity and outgoing flow. The incoming flow is assumed to be at setpoint traffic flux. We do not assume any prior knowledge of the initial conditions and set the initial conditions to be at the setpoint density and velocity. We can see that state estimates converges to the values of plant after $75 \; \rm s$. 

In Fig. 3, the evolution of estimation errors are shown. After $75 \; \rm s$, the estimation errors for density and velocity converge to value less than $1\%$ of the setpoint value. There are still relatively very small estimation errors remain in the domain for two reasons. Our result only guarantees the convergence of estimates in the spatial $L^2$ norm. In addition, there could be nonlinearities of the error system not driven to zero by the linear output injections of the nonlinear boundary observer design.

\section{Data Validation}

In this section, we validate our boundary observer design with Next Generation Simulation (NGSIM) traffic data ~\cite{NGSIM} which provides vehicle trajectories with great details and accuracy. The NGSIM trajectory data set is collected on April 13, 2005 by the Federal Highway Administration's project. The study area is a segment of Intestate 80 located at Emeryville, California. The dataset gathers trajectories of vehicles over a total of 45 minutes during rush hour: 4:00pm - 4:15pm, 5:00pm - 5:15pm, 5:15pm - 5:30pm.

Firstly, we calibrate the nonlinear ARZ model with part of the NGSIM data to obtain calibrated model parameters including  the steady state values, the equilibrium velocity-density function $V(\rho)$ and the relaxation time $\tau$. Then the rest datasets are used to test the observer design for the calibrated ARZ model. The estimation results of traffic states are compared with the NGSIM data.  The boundary data is measured directly from the NGSIM data and traffic states are estimated for the considered domain. The result of reconstructed traffic data and boundary observer estimation of the traffic states are compared.

\subsection{Model calibration with NGSIM data}

\subsubsection{Reconstruction from Data}
We aim to calibrate the ARZ model which is a macroscopic model describing aggregated values. However, the NGSIM data set consists of microscopic measurements. The data was recorded with high-speed cameras for every 0.1 seconds. We need to process NGSIM trajectory data into macroscopic scale so that it can be used to calibrate the ARZ model. 

%Note that here we can also use detector data set from the freeway Performance Measurement System (PeMS) project. This data set provides macroscopic fundamental diagram data, which is obtained by loop detector to record averaged FD data over 30s.

The data was recorded on a $537$-meter long freeway segment with six lanes for a time period of $15$-minutes. Due to insufficient data collection at boundaries of segment, onset and offset of recording, the viable domain we choose to use in calibration and validation is $400$-meter during a time period around $10$-minutes. When we calibrate the parameters in ARZ model and fundamental diagram, we consider the freeway segment as a macroscopic general one-lane problem. That being said, six-lane densities need to be taken into account. 

We will use the Edie's formula \cite{Edie} to calculate aggregated traffic states $\rho(x,t), v(x,t), q(x,t)$ from the trajectory data of vehicles $x(t)$ with a resolution $0.1\rm\; s$. At each time instance, positions of the multiple vehicles are collected. Consider a time-space domain $[0,T] \times [0,L] $, we divide it into $N \times M$ grids 
\begin{align*}
[i \Delta t, (i+1) \Delta t] \;\times\; [j \Delta x, (j+1) \Delta x],
\end{align*} 
where $i \in {1,2,..,N}$ and $j \in  {1,2,..,M}$.
Within each cell, we consider $\rho_{i,j}, q_{i,j}, v_{i,j}$ to be constant. We use the following Edie's formula to map a set of vehicles' traces to speed, flow and density over the space-time grid. For each cells, suppose there are $N_{ij}$ vehicle traces passing through the cell $[i \Delta t, (i+1) \Delta t] \;\times\; [j \Delta x, (j+1) \Delta x],$
\begin{align}
\rho_{i,j} =& \frac{\Sigma_{k=1}^{N_{ij}} t_{k} }{ \Delta x  \Delta t},\\
q_{i,j} =& \frac{\Sigma_{k=1}^{N_{ij}} x_{k} }{ \Delta x  \Delta t},\\
v_{i,j} =& \frac{q_{i,j}}{\rho_{i,j}}.
\end{align}
After obtaining the cell values $\rho_{i,j}, q_{i,j}, v_{i,j}$, they can be later on compared with the observer estimates $\hat \rho_{i,j}, \hat q_{i,j}, \hat v_{i,j}$ with same griding. The number of cells are chosen such that in each cell, there are enough trajectory data. Otherwise, there could be cells that no trajectory has crossed. On the other hand, noises appear if a very fine discretization of grids is chosen. The following simulation is performed in a $ 41 \times 41$ grid. 

\begin{figure*}[t!]
	\centering
	\includegraphics[width=14cm]{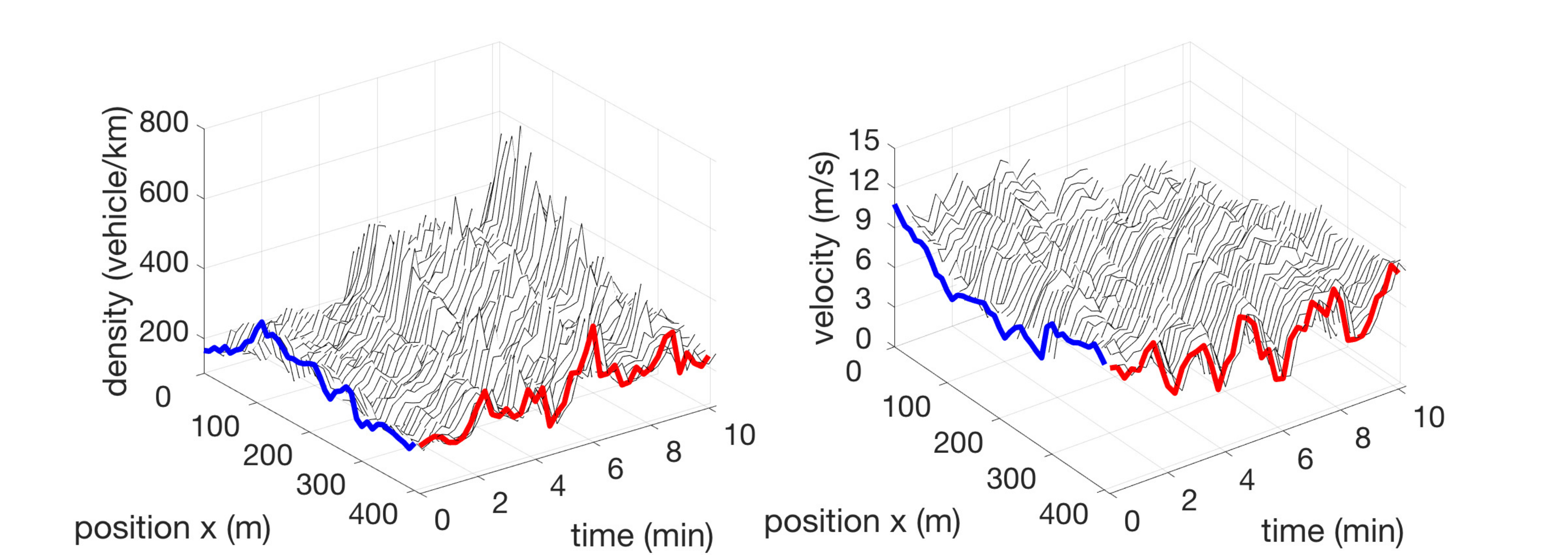}
	\caption{\label{Fig4pm}Density and velocity reconstructed from data of 4:00pm-4:15pm.}
\end{figure*}

\begin{figure*}[t!]
	\centering
	\includegraphics[width=14cm]{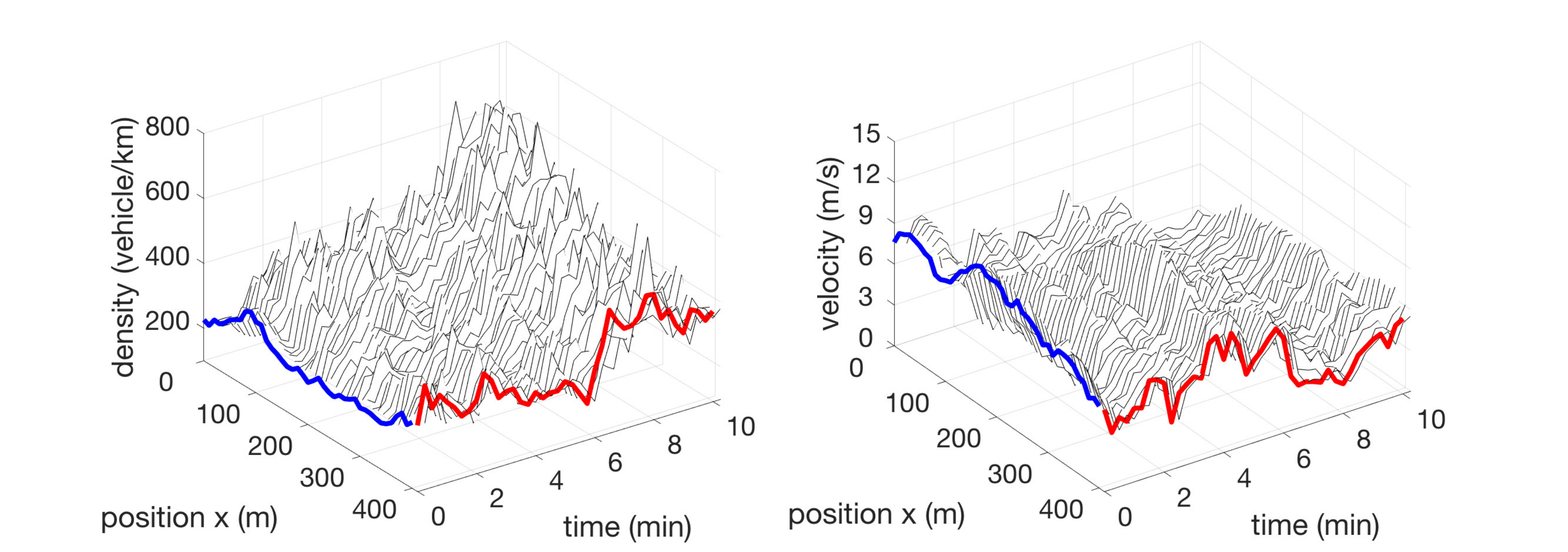}
	\caption{\label{Fig5pm}Density and velocity reconstructed from data of 5:00pm-5:15pm.}
\end{figure*}

We reconstruct the aggregated traffic states from all the three dataset. In Fig.~\ref{Fig4pm} and  Fig.~\ref{Fig5pm}, we show the surface plot of the density and velocity states for the dataset of 4:00pm - 4:15pm and the dataset of 5:00pm - 5:15pm. The initial conditions are highlighted with color red and the boundary conditions at outlet are highlighted with color blue. The congestion forms up as time goes by and propagates from the downstream to upstream. The most congested traffic appears at the inlet where the traffic density is relatively high and velocity is low. 

We are mostly interested in the congested traffic where estimation of the traffic states becomes more relevant. The linearized ARZ model around the uniform reference is analyzed and employed for the observer design. By taking average of traffic aggregated values, we obtain the reference system $\rho^\star$, $v^\star$ and $q^\star$ of each dataset.
 Therefore, the average density, average velocity and average flow of each time period is calculated and shown in the Table~\ref{ave}. 
We observe that among the three data set, the traffic is most congested during 5:15pm - 5:30pm with largest averaged density and smallest velocity. Whether the traffic states are in congested or free regime need to be determined after we introduce the calibrated fundamental diagram.

\begin{table}[t!]
	\centering
	\caption{\label{ave}Averaged aggregate traffic data}
	\begin{tabular}{ p{2cm}|p{1.4cm}|p{1.4cm}|p{1.4cm}}
		\hline
		\bf Data Set & \bf  Density  \rm (veh/km)& \bf  Velocity \rm (km/h) & \bf  Flow  \rm (veh/h)\\ \hline
		4:00 - 4:15pm &  267 & 28.27 & 7548 \\  \hline
		5:00 - 5:15pm & 353 & 20.23  & 7141 \\ \hline  
		5:15 - 5:30pm & 375 & 19.35  & 7256 \\ \hline
	\end{tabular}
\end{table}

\subsubsection{Calibration of model parameters}
For the ARZ model,
\begin{align}
\partial_t \rho + \partial_x( \rho v)=&0, \\
\partial_t v+(v + \rho V'(\rho))\partial_x v=&\frac{V(\rho)-v}{\tau},
\end{align}
the model parameters to be calibrated from the dataset is the equilibrium density-velocity relation $V(\rho)$ and relaxation time $\tau$. The fundamental diagram defined as
\begin{align}
	Q(\rho) = \rho V(\rho)
\end{align}
describing the equilibrium density and flow rate relation is usually obtained by long-term measurements via loop-detectors. The loop-detector data set provides macroscopic density and flow rate data and its recording resolution is $30 \;\rm s$. In the previous section, we use Greenshield's model \eqref{vf} for $V(\rho)$ as a simple choice for the boundary observer design. The Greenshield's fundamental diagram $Q(\rho)$ is given by
\begin{equation}
Q(\rho)=\rho v_f\left(1-\left(\frac{\rho}{\rho_m}\right)^\gamma\right). \label{qr}
\end{equation}
But Greenshield's model cannot accurately represent the fundamental diagram data. The critical density $\rho_c$ satisfies $Q'(\rho)|_{\rho_c} = 0$ and thus segregates the free and congested regimes. The critical density $\rho_c$ of the Greenshield's model $(\gamma =1)$ occurs at $\rho_c = \frac{1}{2} \rho_m$. However, the critical density obtained from empirical traffic data usually shows up at  $\rho_c = \frac{1}{4} \rho_m$. Hence, we need to consider a more realistic functional form for $Q(\rho)$. Here we employ a three-parameter fundamental diagram proposed by \cite{Fan:data}. 

%
%\begin{figure}[h!]
%	\includegraphics[width=14cm]{p1}
%	\centering
%	\caption{Comparison of three-parameter fundamental diagram and Greenshield's model}
%\end{figure} 

\begin{figure*}[t!]
	\centering
	\includegraphics[width=18cm]{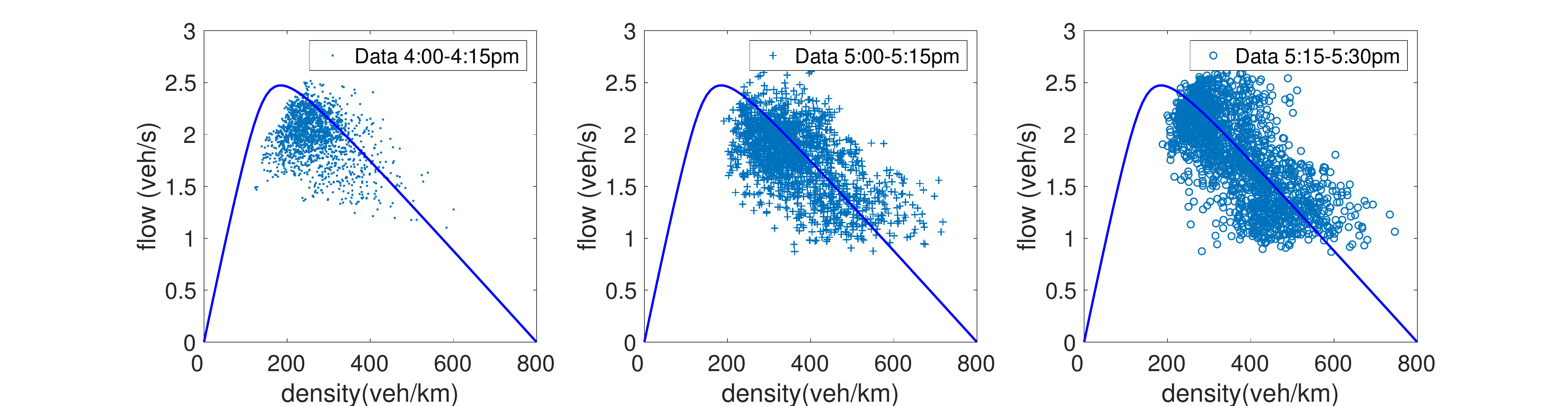}
	\caption{\label{FD} Density and flow from data of 4:00pm-4:15pm, 5:00pm-5:15pm and 5:15pm-5:30pm. }
\end{figure*}
In \cite{Fan:data}, the following three-parameter $(\lambda, p, \alpha)$ fundamental diagram is calibrated with the NGSIM detector data set of the same freeway segement, 
\begin{align}
Q(\rho) = \alpha \left(a + (b-a) \frac{\rho}{\rho_{m}} - \sqrt{1+\lambda^2\left(\frac{\rho}{\rho_{m}}-p\right)^2}\right),
\end{align}
where $a$ and $b$ are denoted by
\begin{align}
	a =& \sqrt{1+(\lambda p)^2},\\
	b =& \sqrt{1+(\lambda (1-p))^2}.
\end{align}
The parameters $(\lambda, p, \alpha)$ do not have physical meaning but represent the shape of the functional form where $\lambda$ represents the roundness, $p$ tunes the critical density, $\alpha$ determines the maximum flow rate. The hyperbolicity $Q''(\rho)<0, V'(\rho)<0$ is guaranteed. 
The three parameters $(\lambda, p, \alpha)$ are determined using Least Square fitting with historical loop detector data.
\begin{figure*}[t!]
	\centering
	\includegraphics[width=15cm]{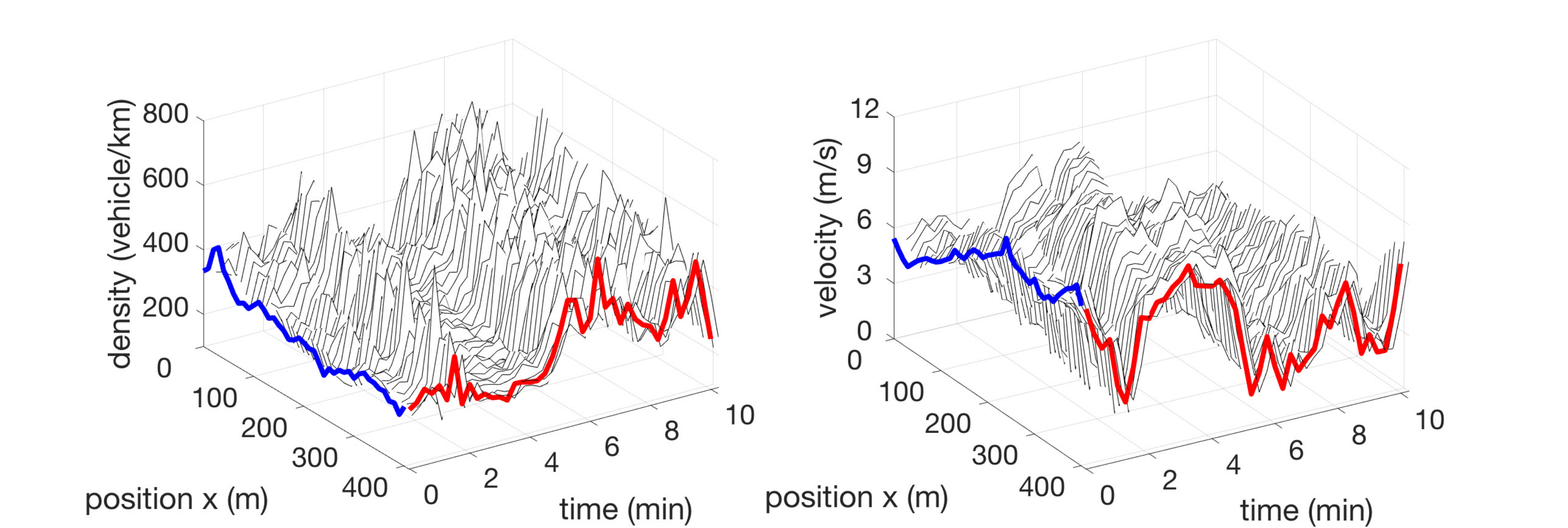}
	\caption{\label{fig:data}Density and velocity reconstructed from the data of 5:15pm-5:30pm.}
	\includegraphics[width=15cm]{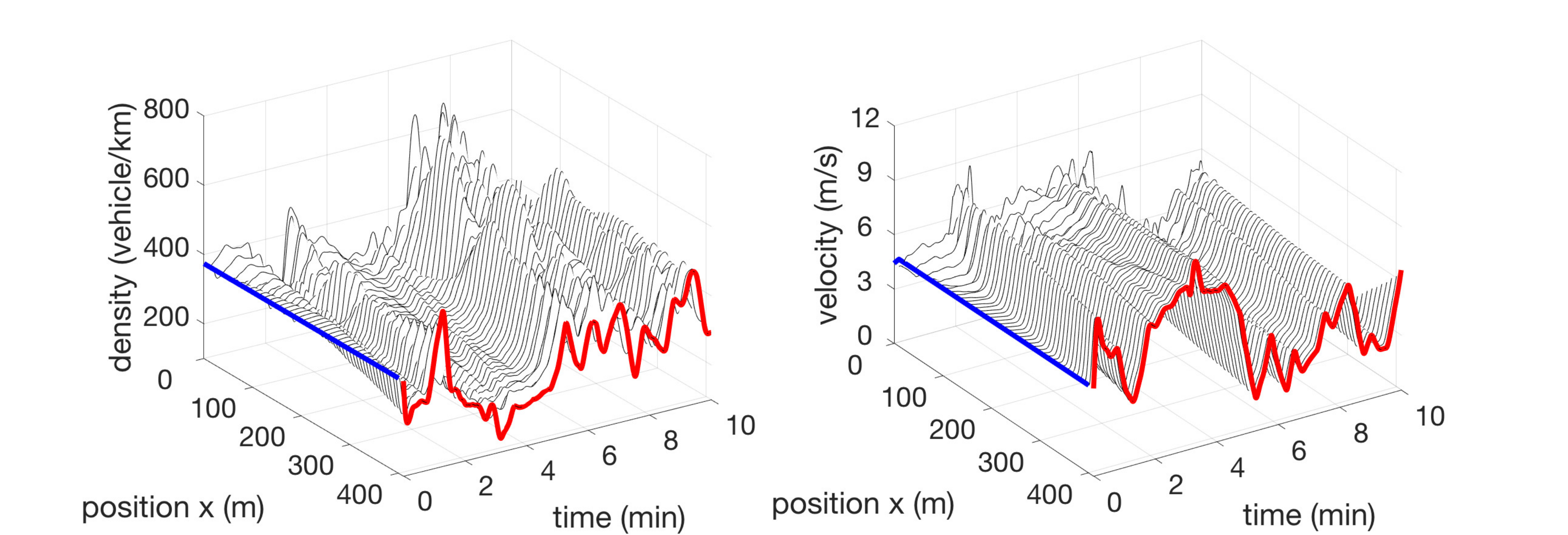}
	\caption{\label{fig:data_estimate}Estimates of density and velocity from the data of 5:15pm-5:30pm.}
\end{figure*}

Due to the lack of data near the maximum density, the value of $\rho_m$ is prescribed according to the following equation
\begin{align}
\rho_m = \frac{\rm number\; of\; lanes}{\rm typical\; vehicle\; length \times safety\; distance \;factor}.
\end{align}
The freeway segment in the dataset consists of $6$ lanes and we consider the typical vehicle length to be $5$-meter and the safety distance factor is $50 \%$ of vehicle length. Therefore, we have $\rho_m$ for all lanes in our simulation
\begin{align}
\rho_m = 800 \; {\rm veh/km}.
\end{align}

The calibrated fundamental diagram is plotted in Fig.~\ref{FD}. The traffic density and flow rate of the three dataset are plotted on the calibrated fundamental diagram. We can see that 4:00pm-4:15pm are in the transition region where the data points are partially in the free regime and partially in the congested regime. The traffic data of 5:00pm-5:15pm and 5:15pm-5:30pm are scattered in the congested regime of the fundamental diagram. 

%I also need to calibrate $\tau$ then do the error analysis. Then with correct parameters, First test model validity and then control result. 

%The boundary observer can be tested firstly but how to implement the metering would be a question to ask. 
With the calibrated fundamental diagram $V(\rho)$, we choose the relaxation time $\tau$ from a range from $10s$ to $100s$ and calibrate it with the dataset of 5:00pm-5:15pm. The optimal relaxation time is $\tau = 30s$ where the total error between the calibrated model and data is the lowest. In the next step, we use the calibrated fundamental diagram $V(\rho)$ and the relaxation time $\tau$ to construct the boundary observer. 
%
%\begin{figure*}[h!]
%	\centering
%	\includegraphics[width=18cm]{model}
%	\caption{ARZ model prediction with calibrated parameters from NGSIM data}
%	\includegraphics[width=18cm]{data}
%	\caption{NGSIM traffic data}
%\end{figure*}

\subsection{Simulation for the nonlinear observer with calibrated parameters}

We use the data of 5:15pm-5:30pm to test the boundary observer design. The reference system $(\rho^\star, v^\star, q^\star)$ is obtained from Table~\ref{ave}. Along with the calibrated parameters $V(\rho)$ and $\tau$, the nonlinear observer is constructed with a copy of the nonlinear ARZ model with the output injection gains that drive the estimation errors to zero. The numerical solution of the nonlinear PDEs are approximated with the Lax-Wendroff method. The boundary data is implemented with the ghost cell. The ARZ model collects the boundary values based both on flux of the computational domain and the boundary data of the ghost cells. Using the boundary measurements of the inlet and outlet of the freeway segment, the state estimation $(\hat \rho(x,t), \hat v(x,t))$ is generated without the knowledge of the initial condition. In Fig.~\ref{fig:data}, $(\rho(x,t),  v(x,t))$ is obtained from the reconstruction of the data set of 5:15pm-5:30pm. In Fig.~\ref{fig:data_estimate}, it shows the evolution of the state estimates $(\hat \rho(x,t), \hat v(x,t))$. The initial condition, highlighted with color blue, is assumed to be the uniform reference system $(\rho^\star, v^\star, q^\star)$ which represents the averaged values of the dataset. The boundary conditions at outlet are highlighted with right color which gives the output injections in the observer. We notice that when density value is higher than $600 \;$veh/km at inlet around $7\; \rm\min$, the estimation result is not satisfying at inlet. This could be related to the  ARZ model's inaccuracy in predicting traffic states near maximum density since non-unique maximum densities exist for the ARZ model. 

For the error analysis of the observer estimation, the estimation errors are considered in the $L^2$-norm, defined as
\begin{align}
	E_{\rho}(t) =& \left[\frac{1}{L}\int_{0}^{L}\left(\frac{\rho(x,t) - \hat\rho(x,t)}{\rho^\star}\right)^2 dx\right]^{1/2},\\
	E_{v}(t) =& \left[\frac{1}{L}\int_{0}^{L}\left(\frac{v(x,t) - \hat v(x,t)}{v^\star}\right)^2 dx\right]^{1/2},	
\end{align}
where $\rho^\star$ and $v^\star$ are the averaged state values of the data. We choose the $L^2$ of the estimation errors and average it over space. 
The convergence of the local stability in the $L^2$-sense for estimation errors to zero is guaranteed in Theorem 2. In addition, the spatial averaged errors can remove the influence of noises and outliers of the traffic data. 

\begin{figure}[t!]
	\centering
	\includegraphics[width=9.6cm]{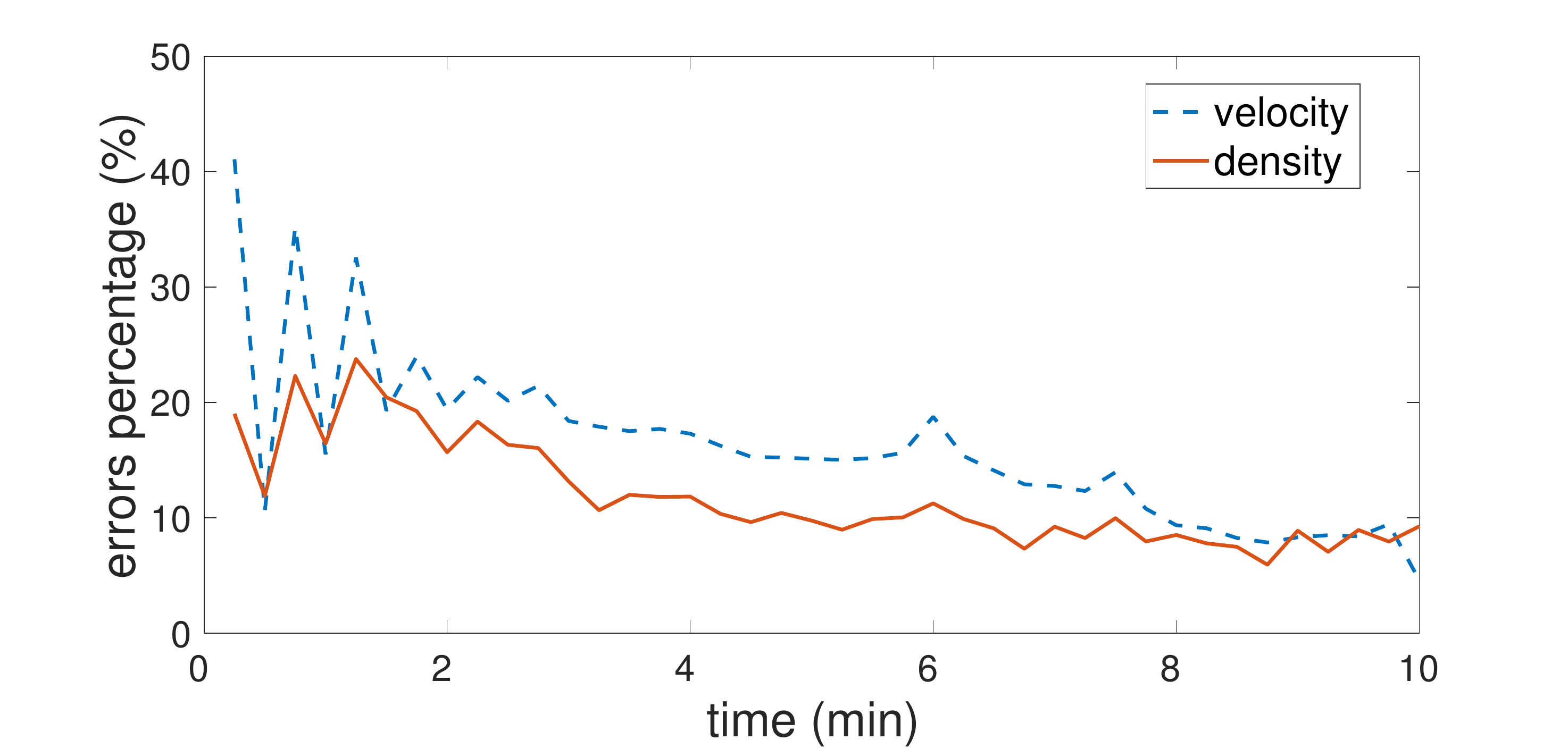}
	\caption{\label{fig:error}Density and velocity estimation errors for the data of 5:15pm-5:30pm.}
\end{figure}

The temporal evolution of the space-averaged errors of density and velocity estimates in the $L^2$-sense is shown in Fig.~\ref{fig:error}. It reveals that at the initial time, density and velocity estimation errors start from $20\%$ and
$40\%$ respectively. The finite convergence time is around $t_f = 3 \; $min. The estimation errors in the end converge at $10\%$. The linearization of output injections design, the data noise, the reconstruction errors and the numerical approximation errors could contribute to the remaining spatial averaged errors between the estimation and NGSIM traffic data after the convergence time.

%simulation \\
%writing theorem \\
%compare model predicted result, linearized observer, nonlinear observer. maybe kalman fitler \\

%I210 including both congested 
%
%\begin{itemize}
%	\item $V(\rho)$ can be chosen arbitrary
%	\item nonlinear observer design has not been considered before
%	\item numerical validation with real data has never been performed for nonlinear observer
%\end{itemize}

\section{Conclusion}

In conclusion, we develop a nonlinear boundary observer for the second-order nonlinear hyperbolic PDEs, estimating traffic states of ARZ model and then validate the design with traffic field data. Analysis of the linearized ARZ model leads our main focus to the congested regime where stop-and-go happens. Using spatial transformation and PDE backstepping method, we construct a boundary observer with a copy of the nonlinear plant and output injection of measurement errors so that the exponential stability of estimation errors in the $L^2$ norm and finite-time convergence to zero are guaranteed. Simulations are performed for traffic estimation on a stretch of freeway. The nonlinear observer is tested with a calibrated ARZ model obtained from the NGSIM data. 

For future work, observer design may be considered for a generalized ARZ model proposed by \cite{Seibold3} to address the non-unique maximum density associated with ARZ model. The estimation accuracy in predicting the heterogeneous behaviors of drivers and spread of data for the congested regime could be improved. On the other hand, defining the fundamental diagram requires the calibration with the historical data. This assumption of using the historical data to determine model parameters may not hold when traffic becomes unpredictable in case of accidents. It is practically preferable if the model parameters could be estimated real-time. Therefore, it is of authors' interest to consider adaptive observer design for this problem.

                                                                         % in the appendices.
\end{document}